\documentclass[twoside,leqno,twocolumn]{article}  
\usepackage{ltexpprt}

\usepackage{doi,url}
\hypersetup{
colorlinks,
citecolor=blue,
linkcolor=black,
urlcolor=blue}

\usepackage{amsmath}
\usepackage{amssymb} % for \begin{bmatrix}, \begin{smallmatrix}, \begin{align} etc.
\usepackage{amsfonts} % for \mathbb{} etc. `blackboard' for uppercase only.
\usepackage{mathrsfs}

\usepackage{scalefnt}
\usepackage{graphicx}
\usepackage{subfig}

\usepackage{algorithm,algorithmic}

\begin{document}

\title{\Large Preconditioned Continuation Model Predictive Control}
\author{Andrew Knyazev\thanks{Mitsubishi Electric Research Labs (MERL)
201 Broadway, 8th floor, Cambridge, MA 02139, USA,
\href{mailto:knyazev@merl.com}{knyazev@merl.com}, 
\url{http://www.merl.com/people/knyazev}.} \\
\and
Yuta Fujii\thanks{This work was performed while the author was affiliated
with Advanced Technology R\&D Center, Mitsubishi Electric Corporation, 8-1-1,
Tsukaguchi-honmachi, Amagasaki City, 661-8661, Japan,
http://www.mitsubishielectric.com/company/rd/research/labs/
advanced\_technology/.
}
\and
Alexander Malyshev\thanks{Mitsubishi Electric Research Labs (MERL)
201 Broadway, 8th floor, Cambridge, MA 02139, USA,
\href{mailto:malyshev@merl.com}{malyshev@merl.com},  
\url{http://www.merl.com/people/malyshev}.}}
\date{}

\maketitle
 
%\pagenumbering{arabic}
%\setcounter{page}{1}%Leave this line commented out.

\begin{abstract} \small\baselineskip=9pt
Model predictive control (MPC) anticipates future events to take appropriate
control actions. Nonlinear MPC (NMPC) describes systems with nonlinear models
and/or constraints. A Continuation/GMRES Method for NMPC, suggested by T. Ohtsuka
in 2004, uses the GMRES iterative algorithm to solve a forward difference approximation
$Ax=b$ of the Continuation NMPC (CNMPC) equations on every time step. The coefficient matrix
$A$ of the linear system is often ill-conditioned, resulting in poor GMRES convergence,
slowing down the on-line computation of the control by CNMPC, and
reducing control quality. We adopt CNMPC for challenging minimum-time problems, and
improve performance by introducing efficient preconditioning, utilizing parallel computing, 
and substituting MINRES for GMRES. 
\end{abstract}

\section{Introduction}

Model predictive control (MPC) is used in many applications to control complex dynamical systems. Examples of such systems include production lines, car engines, robots, other numerically controlled
machining, and power generators. The MPC is based on optimization of the operation of the system
over a future finite time-horizon, subject to constraints, and implementing the control only over the current time step.
 
Model predictive controllers rely on dynamic models of the process, most often linear empirical models, in which case the MPC is linear. Nonlinear MPC (NMPC), which describes systems with nonlinear
models and constraints, is often more realistic, compared to the linear MPC, but computationally more difficult. Similar to the linear MPC, the NMPC requires solving optimal control problems on a finite prediction horizon, generally not convex, which poses computational challenges.

Numerical solution of the NMPC optimal control problems may be based on Newton-type optimization schemes. Exact Newton-type optimization schemes require an analytic expression of a corresponding Jacobian matrix, which is rarely available in practice and is commonly replaced with a forward difference (FD) approximation; see, e.g., \cite{K95}. Such approximate Newton-type optimization schemes utilize the FD approximation of the original nonlinear equation during every time step.
An efficient variant of the approximate Newton-type optimization can be performed
by a Continuation NMPC (CNMPC) numerical method proposed by T.~Ohtsuka in \cite{O04},
where each step of the algorithm requires solving a system of linear equations performed
by the GMRES iterative method \cite{SS86}.

Our contributions presented below are two-fold. We describe an extension of CNMPC with a terminal constraint, suitable to solve 
minimum-time optimal control problems, and with an optimization parameter. 
We~investigate preconditioning
for GMRES in the context of the NMPC problems and using the MINRES iteration
\cite{PS74} instead of GMRES.
MINRES provides overall faster implementation,
compared to GMRES without restarts, of our approach in cases, where many iterations are required.  
Our numerical simulations show that
the preconditioning can considerably improve the quality of controllers
with marginal extra computational time, which can be reduced or eliminated by employing a parallel processing
for the preconditioner setup. 

The rest of the paper is organized as follows. In Section \ref{s2}, we formulate CNMPC of Ohtsuka, extended to having a terminal constraint and a parameter. Section~\ref{s3} describes the original algorithm of Ohtsuka, where the FD linear system is solved using GMRES, and then introduces MINRES as an alternative to GMRES, discusses preconditioning  for GMRES and MINRES, and suggests specific algorithms of constructing the preconditioner and using it to accelerate convergence of iterations. In Section \ref{s4}, we give a detailed description of a test minimum-time optimal control problem, defining  a quickest  arrival of the system to a given destination, with inequality constraints on the system control, and its CNMPC formulation.  Section \ref{s5} presents our results of numerical experiments solving the test problem, demonstrating advantages of the proposed approaches.

\section{Finite horizon optimization by CNMPC}\label{s2}

As a specific example of a mathematical formalism of NMPC, we consider an extended
version of the control problem considered by T. Ohtsuka \cite{O04} as follows,
\[
\min_{u,p} J,
\]
\[
J = \phi(t+T,x(t+T),p)+\int_t^{t+T}L(t',x(t'),u(t'),p)dt'
\]
subject to
\begin{equation}\label{e1}
\dot{x} = \frac{dx}{dt'}=f(t',x(t'),u(t'),p),
\end{equation}
\begin{equation}\label{e2}
C(t',x(t'),u(t'),p) = 0,
\end{equation}
\begin{equation}\label{e3}
\psi(t+T,x(t+T),p) = 0.
\end{equation}

Here, $x=x(t)$ denotes the vector of the state of the dynamic system, also serving
as an initial state for the optimal control problem over the horizon.
The vector $u=u(t)$ is the control vector, serving as an input to control the system.
The scalar function $J$ describes a performance cost to be minimized,
which includes a terminal cost (the first term
in the sum) and a cost over the finite horizon (the second term in the
sum). Equation (\ref{e1}) is the system dynamic model that may be nonlinear in $x$
and/or $u$. Equation (\ref{e2}) describes the equality constraints for the state $x$
and the control~$u$. 
The horizon time length $T$ may in principle also depend on $t$, e.g.,\ for time-optimal control problems.
In this case,
the original problem can be converted into a fixed horizon problem by letting
$T(t)=1\cdot t_f$, where $t_f$ is an additional parameter to be included in $p$
and determined in MPC. Substituting $t+\tau t_f$ for the time $t'$, we arrive at
a problem with the normalized time scale $\tau$ and fixed horizon $[t,t+1]$.
Such a conversion is applied to the test problem in Section \ref{s4}.

Compared to \cite{O04}, one extra constraint (\ref{e3}),
described by the terminal constraint function $\psi$, and an extra parameter vector $p$ are being added 
to the problem formulation, allowing one to extend CNMPC to a 
wide range of optimal control and design problems.

The NMPC optimal control problem is solved by a variational approach. Its discrete
counterpart is solved by the traditional Lagrange method of undetermined multipliers.
We denote the costate vector by $\lambda$ and the Lagrange multiplier vector associated
with the equality constraint (\ref{e2}) by $\mu$. The terminal constraint (\ref{e3})
is relaxed by introducing the Lagrange multiplier $\nu$.
The so-called Hamiltonian function, as defined in control theory, is
\begin{eqnarray*}
\lefteqn{H(t,x,\lambda,u,\mu,p) = L(t,x,u,p)}\hspace*{6em}\\
&&{}+\lambda^Tf(t,x,u,p)+\mu^TC(t,x,u,p).
\end{eqnarray*}

To discretize the continuous formulation of the optimal control problem stated above,
we introduce a uniform horizon time grid by dividing the horizon $[t,t+T]$ into $N$
time steps of size $\Delta\tau$ and replace the time-continuous vector functions $x(\tau)$
and $u(\tau)$ by their indexed values $x_i$ and $u_i$ at the grid points. Thus, $N$ is
a number of artificial time steps for the optimal control problem over the horizon.
The integral in the performance cost $J$ over the time horizon is approximated
by a simple quadrature rule. The time derivative of the state vector is approximated
by the forward difference formula. Then the discretized optimal control problem
appears as follows,
\[
\min_{u_i,p} J,
\]
\[
J = \phi(\tau_N,x_N,p) + \sum_{i=0}^{N-1}L(\tau_i,x_i,u_i,p)\Delta\tau,
\]
subject to
\begin{equation*}%\label{e4}
\qquad x_{i+1} = x_i + f(\tau_i,x_i,u_i,p)\Delta\tau,\quad i = 0,1,\ldots,N-1,
\end{equation*}
\begin{equation*}%\label{e5}
C(\tau_i,x_i,u_i,p) = 0,\quad  i = 0,1,\ldots,N-1,
\end{equation*}
\begin{equation*}%\label{e6}
\psi(\tau_N,x_N,p) = 0.
\end{equation*}
We note that we have so far discretized the NMPC optimal control problem
only in the horizon time. We will discretize the system time $t$ later using the uniform
time step size $\Delta t$, i.e. discretization in the horizon time may be different
from the time discretization of the system.

The necessary optimality conditions for the discretized horizon problem
are obtained using the discrete Lagrangian function
\begin{eqnarray*}
\lefteqn{\mathscr{L}(X,U)=\phi(\tau_N,x_N,p)+\sum_{i=0}^{N-1}L(\tau_i,x_i,u_i,p)\Delta\tau}\\
&&{}+\lambda_0^T[x(t)-x_0]\\
&&+\sum_{i=0}^{N-1}\lambda_{i+1}^T[x_i-x_{i+1}
+f(\tau_i,x_i,u_i,p)\Delta\tau]\\
&&+\sum_{i=0}^{N-1}\mu_i^TC(\tau_i,x_i,u_i,p)\Delta\tau+\nu^T\psi(\tau_N,x_N,p),
\end{eqnarray*}
where $X = [x_i\; \lambda_i]^T$ and $U = [u_i\; \mu_i\; \nu\; p]^T$.
Namely, the necessary optimality conditions coincide with the stationarity conditions
\[
\frac{\partial \mathscr{L}^T}{\partial X}(X,U)=0 \mbox{ and }
\frac{\partial \mathscr{L}^T}{\partial U}(X,U)=0.
\]
For example, the derivative with respect to $u_i$, which is
$\partial \mathscr{L}^T/\partial u_i=0$, yields the following equation:
\begin{eqnarray*}
\lefteqn{\frac{\partial L}{\partial u_i}(\tau_i,x_i,u_i,p)\Delta\tau+\lambda_{i+1}^T
\frac{\partial f}{\partial u_i}(\tau_i,x_i,u_i,p)\Delta\tau}\hspace*{6.2em}\\
&&+\mu_i^T\frac{\partial C}{\partial u_i}(\tau_i,x_i,u_i,p)\Delta\tau=0.
\end{eqnarray*}
Using the Hamiltonian function, it can be shortened to
\[
\frac{\partial H}{\partial u_i}(\tau_i,x_i,\lambda_{i+1},u_i,\mu_i,p)\Delta\tau=0.
\]
Taking the derivative with respect to $\mu_i$, which is
${\partial \mathscr{L}^T}/{\partial\mu_i}=0$, we obtain the following equation,
which also involves the factor $\Delta\tau$,
\[
C(\tau_i,x_i,u_i,p)\Delta\tau = 0.
\]

Now we proceed to the construction of a vector function $F(U,x,t)$,
which is used to formulate the full set of necessary optimality conditions.
The vector function $U=U(t)$ combines the control input $u$, the Lagrange multiplier $\mu$,
the Lagrange multiplier $\nu$, and the parameter $p$, all in one vector, as follows, 
\[
U(t)=[u_0^T,\ldots,u_{N-1}^T,\mu_0^T,\ldots,\mu_{N-1}^T,\nu^T,p^T]^T. 
\]
The vector argument $x$ in the function $F(U,x,t)$ denotes the current measured
state vector, which serves as the initial vector $x_0$ in the following algorithm,
defining an evaluation of  $F(U,x,t)$. 
\begin{enumerate}
\item Starting with the current measured state $x_0$, compute $x_i$, $i=1,2\ldots,N$,
by the forward recursion
\[
x_{i+1} = x_i + f(\tau_i,x_i,u_i,p)\Delta\tau,\, i=0,\ldots,N-1.
\]
Then starting with the value
\[
\lambda_N=\frac{\partial\phi^T}{\partial x}(\tau_N,x_N,p)+
 \frac{\partial\psi^T}{\partial x}(\tau_N,x_N,p)\nu
\]
compute the costate $\lambda_i$, $i=N\!-\!1,\ldots,0$, by the backward recursion
\[
\lambda_i=\lambda_{i+1}+\frac{\partial H^T}{\partial x}
(\tau_i,x_i,\lambda_{i+1},u_i,\mu_i,p)\Delta\tau.
\]
\item Calculate the vector function $F[U,x,t]$, using the just obtained $x_i$ and $\lambda_i$, 
$i=0,1\ldots,N$,
as follows,
\begin{eqnarray*}
\lefteqn{F[U,x,t]}\\
&&\hspace*{-2em}=\left[\begin{array}{c}\begin{array}{c}
\frac{\partial H^T}{\partial u}(\tau_0,x_0,\lambda_{1},u_0,\mu_0,p)\Delta\tau\\
\vdots\\\frac{\partial H^T}{\partial u}(\tau_i,x_i,\lambda_{i+1},u_i,\mu_i,p)\Delta\tau\\
\vdots\\\frac{\partial H^T}{\partial u}(\tau_{N-1},x_{N-1},\lambda_{N},u_{N-1},
\mu_{N-1},p)\Delta\tau\end{array}\\\;\\
\begin{array}{c}C(\tau_0,x_0,u_0,p)\Delta\tau\\
\vdots\\C(\tau_i,x_i,u_i,p)\Delta\tau\\\vdots\\
C(\tau_{N-1},x_{N-1},u_{N-1},p)\Delta\tau\end{array}\\\;\\
\psi(\tau_N,x_N,p)\\[2ex]
\begin{array}{c}\frac{\partial\phi^T}{\partial p}(\tau_N,x_N,p)+
\frac{\partial\psi^T}{\partial p}(\tau_N,x_N,p)\nu\\
+\sum_{i=0}^{N-1}\frac{\partial H^T}{\partial p}(\tau_i,x_i,
\lambda_{i+1},u_i,\mu_i,p)\Delta\tau\end{array}
\end{array}\right].
\end{eqnarray*}
\end{enumerate}

The optimality condition is the nonlinear equation
\begin{equation}\label{e7}
 F[U(t),x(t),t]=0
\end{equation}
with respect to the unknown $U(t)$, which needs
to be solved numerically by a computer processor at each time step of NMPC
in real time on the controller board. This is the most difficult and
challenging part of implementation of NMPC.  
At the initial time $t=t_0$, we need to approximately solve \eqref{e7} directly. 

Let us denote the step size of the system time discretization by $\Delta t$, 
assume that $U(t-\Delta t)$ is already available at the time $t$, and set 
$\Delta U=U(t)-U(t-\Delta t).$
For a small scalar $h>0$, which may be different from the system time step
$\Delta t$ and from the horizon time step~$\Delta \tau$, we introduce the operator
\begin{eqnarray}\label{e9}
\lefteqn{a(V)=(F[U(t-\Delta t)+hV,x(t),t]}\hspace*{4em}\\
&&{}-F[U(t-\Delta t),x(t),t])/h.\nonumber
\end{eqnarray}
Then equation (\ref{e7}) is equivalent to the equation 
\begin{equation*}%\label{e8}
h a(\Delta U/h)=b, \text{ where } b=-F[U(t-\Delta t),x(t),t].
\end{equation*}

Let us denote the $j$-th column of the $m\times m$ identity matrix by $e_j$,
where $m$ is the dimension of the vector $U$, and construct an $m\times m$ matrix $A$
with the columns $Ae_j$, $j=1,\ldots,m$, defined by the formula
%\begin{eqnarray}\label{e10}
 %\lefteqn{\qquad Ae_j=\left(F[U(t-\Delta t)+he_j,x(t),t]\right.}\hspace*{5em}\\
 %&&\left.{}-F[U(t-\Delta t),x(t),t]\right)/h.\nonumber
%\end{eqnarray}
\begin{equation}\label{e10}
Ae_j=a(e_j).
\end{equation}
The matrix $A$ approximates the symmetric Jacobian matrix $F_U[U(t-\Delta t),x(t),t]$
so that $a(V)=AV+O(h)$.

It is important to realize that the operator $a(\cdot)$ in~\eqref{e10}
may be nonlinear. In particular, this explains why our algorithms of explicitly
computing $A$ for the purpose of a preconditioner setup may result in
a non-symmetric matrix $A$. 
Numerical stability of computations may be improved by enforcing the symmetry,
by substituting $(A+A^T)/2$ for $A$.
The deviation from the symmetry gets smaller
with a sampling period $h$, which we are free to choose independently of $\Delta t$ and $\Delta \tau$. 

A key limitation
in the choice of $h$ comes from the fact that the cancellation error starts
picking up in the finite difference evaluation in the operator $a(V)$
due to inexact arithmetic of the controller processor.
This is an unavoidable side effect of using the finite difference
approximation of the derivative. A recommended lower bound for the value
of $h$ can for example be $10^{-8}$ in the double precision arithmetic,
but the optimal value also depends on the function $F[U,x,t]$.

Given the formulas for computing the vector function $F[U,x,t]$, nonlinear
equation (\ref{e7}) must be solved at the points of the grid $t_i=t_0+i\Delta t$,
$i=0,1,\ldots$.

At the initial state $x_0=x(t_0)$, we find an approximate  solution $U_0$ to the equation $F[U_0,x_0,t_0]=0$ by a suitable optimization procedure. The dimension of the vector $u(t)$
is denoted by $n_u$. Since
\[
U(t)=[u_0^T,\ldots,u_{N-1}^T,\mu_0^T,\ldots,\mu_{N-1}^T,\nu^T,p^T]^T, 
\]
the first block entry of $U_0$, formed from the first $n_u$ elements of $U_0$,
is taken as the control $u_0$ at the state $x_0$.
The next state $x_1=x(t_1)$ is either measured by a sensor or computed
by the formula $x_1=x_0+\Delta tf(t_0,x_0,u_0)$; cf. \eqref{e1}. Now we start the recursion as follows.

At the time $t_i$, where $i>0$, we arrive with the state $x_i$ and the vector $U_{i-1}$. The operator
\[
a_i(V)=\left(F[U_{i-1}+hV,x_i,t_i]-F[U_{i-1},x_i,t_i]\right)/h,
\]
defined by \eqref{e9}, determines an $m\times m$ matrix $A_i$ with the columns
\[
A_ie_j=a_i(e_j),\, j=1,\ldots,m,
\]
as in (\ref{e10}).
At the current time $t_i$, our goal is to solve the following equation
\begin{equation}\label{e11}
 ha_i(\Delta U_i/h)=b_i, \text{where }
%\end{equation}
%with the right-hand side
%\begin{equation}\label{e12}
 b_i=-F[U_{i-1},x_i,t_i].
\end{equation}

Then we set $U_i=U_{i-1}+\Delta U_i$ and choose the first $n_u$ components of $U_i$ as
the control $u_i$. The next state $x_{i+1}=x(t_{i+1})$ either comes from a sensor, estimated, 
or computed by the formula $x_{i+1}=x_i+\Delta tf(t_i,x_i,u_i)$.

Having the basic setup of CNMPC now described, leading to equation \eqref{e11}, 
next we discuss numerical solution of \eqref{e11}. 
Let us highlight that equation \eqref{e11} is never solved exactly in practice, 
thus, a choice of an algorithm may greatly affect not only the performance 
of the controller, but also the computed control as well.

\section{Algorithms}\label{s3}

A direct way to solve \eqref{e11} approximately is generating the matrix $A_i$
and then solving the system of linear equations $A_i\Delta U_i=b_i$ by,
e.g.,\ the Gaussian elimination.

Another way is solving (\ref{e11}) by a suitable Krylov subspace iteration,
e.g.,\ by GMRES \cite{SS86} or MINRES~\cite{PS74} methods, where we do not need to
generate the matrix~$A_i$ explicitly. Namely, we simply use the operator
$a_i(V)$ instead of computing the matrix-vector product $A_iV$,
for arbitrary vectors $V$; cf.,~\cite{K95,KK04}.
In his seminal paper~\cite{O04}, T.~Ohtsuka uses the GMRES iteration.

A typical implementation of the preconditioned GMRES without restarts is given by Algorithm \ref{a1}, where
$Tr$ denotes an action of a precontioner $T$ on a vector $r$, as explained below.
The unpreconditioned GMRES, as in~\cite{O04}, simply uses $z=r$. 
We denote by $H_{i_1:i_2,j_1:j_2}$ the submatrix of $H$ with the entries $H_{ij}$
such that $i_1\leq i\leq i_2$ and $j_1\leq j\leq j_2$. 

\begin{algorithm}
\caption{Preconditioned GMRES without restarts}
\begin{algorithmic}[1]\label{a1}
\REQUIRE $a(v)$, $b$, $x_0$, $k_{\max}$, $T$
\ENSURE Solution $x$ of $a(x)=b$
\STATE $r=b-a(x_0)$, $z=Tr$, $\beta=\|z\|_2$, $v_1=z/\beta$
\FOR{$k=1,\ldots,k_{\max}$}
\STATE $r=a(v_k)$, $z=Tr$
\STATE $H_{1:k,k}=[v_1,\ldots,v_k]^Tz$
\STATE $z=z-[v_1,\ldots,v_k]H_{1:k,k}$
\STATE $H_{k+1,k}=\|z\|_2$
\STATE $v_{k+1}=z/\|z\|_2$
\ENDFOR
\STATE $y=\mbox{arg min}_y\|H_{1:k_{\max}+1,1:k_{\max}}y-[\beta,0,\dots,0]^T\|_2$
\STATE $x=x_0+[v_1,\ldots,v_{k_{\max}}]y$
\end{algorithmic}
\end{algorithm}

We emphasize that the operator $a_i(\cdot)$ may be nonlinear, but approximates
the symmetric Jacobian matrix $F_U[U_{i-1},x_i,t_i]$. This implies a slight deviation from the symmetry property
$V_2^Ta_i(V_1)=(a_i(V_2))^TV_1$ for arbitrary vectors $V_1$ and $V_2$.
We assume that the deviation is small and propose applying the MINRES iteration
to solve equation (\ref{e11}). 
%The conjugate gradient method is not suitable, since
%the operator is generally nondefinite, i.e., the values $V^Ta_i(V)$ may have opposite signs
%for some vectors $V$.

When the operator $a_i(\cdot)$ is linear and symmetric, the projected $(k_{\max}+1)\times k_{\max}$
matrix $H$, constructed by GMRES without preconditioning, is tridiagonal.
The~MINRES method is then a special
variant of GMRES, which makes use of the tridiagonal structure. The table below,
adopted from \cite{CS14}, gives a comparison of 
computational complexities of 
MINRES and GMRES without preconditioning 
for solution of a
linear system $Ax=b$ with a symmetric $m\times m$ matrix $A$ in terms of memory
storage required by working vectors in the solvers and the number of floating-point
operations. By $t_P$ we denote the work needed for evaluating $a_i(V)$.

\vspace{1ex}\hspace{-1.5em}\begin{tabular}{|l|c|c|}
\hline
 Solver&Storage&Work per iteration\\ \hline
 MINRES&$7m$&$t_P+9m$\\ \hline
 GMRES&$(k_{\max}+2)m$&$t_P+(k_{\max}+3)m+\frac{m}{k_{\max}}$\\[.5ex]
 \hline
\end{tabular}
\vspace{1ex}

If the matrix $A_i$ gets ill-conditioned, the convergence of GMRES or MINRES
may stagnate. The convergence can be improved by preconditioning.
A matrix $T_i$ that approximates the matrix $A_i^{-1}$ and such that
computing the product $T_ir$ for an arbitrary vector $r$ is relatively easy,
is referred to as a preconditioner. The preconditioning for the system of linear
equations $Ax=b$ with the preconditioner $T$ formally replaces the original system $Ax=b$
with the equivalent preconditioned linear system $TAx=Tb$. If the condition number
$\kappa(TA)=\|TA\|\|A^{-1}T^{-1}\|$ of the matrix $TA$ is small, convergence of
iterative solvers for the preconditioned system can be fast. However, 
the convergence of the preconditioned GMRES, in contrast to that of the preconditioned MINRES 
with a symmetric positive definite preconditioner, is not necessarily 
determined by the condition number $\kappa(TA)$.
Results on convergence of GMRES in a nonlinear case can be
found in \cite{BM01}.

When the approximate solution $x_{k_{\max}}$ computed by GMRES after $k_{\max}$
iterations is not accurate enough, it is very common to restart GMRES with
$x_0$ equal to $x_{k_{\max}}$ instead of increasing the maximum number of
iterations $k_{\max}$. Practical implementations of GMRES perform 
restarts. Restarts allow to cap the GMRES memory use to $k_{\max}+2$ vectors,
but may significantly slow down the convergence. In our tests, we apply GMRES without restarts
for simplicity of presentation. 

To setup the preconditioner, 
the matrix $A_i$ is computed at some time $t_i$ and then its LU factorization $A_i=LU$
is computed, where $L$ is a lower- and $U$ is an upper-triangular matrix. 
The product $Tr$ is mathematically given by $Tr=U^{-1}(L^{-1}r)$, 
but is computed by back-substitution, which is much cheaper than the computation of the inverses of $L$ and $U$.
The same preconditioner $T$ is used in a number of subsequent grid points
starting from $t_i$. 
The computation of the matrix $A_i$ requires $m$ evaluations $a_i(e_j)$, see \eqref{e10},
that can be efficiently implemented in parallel. 

The symmetry of the preconditioner $T$ can be used to reduce the memory storage
and processor work; see, e.g.,\ \cite{BGL05}. For example, the factorization $T=LDL^T$, see e.g.\ \cite{GVL13}, 
instead of the LU factorization allows us using only half of memory.
The anti-triangular factorization from~\cite{MVD13} may also reduce both the memory
requirements and work in preconditioning. 

MINRES requires symmetric positive definite
preconditioners such as in~\cite{VK13}. In our MINRES simulations, 
although not reported in Section \ref{s5} in details, we use the preconditioned MINRES-QLP method from \cite{CS14}.

\section{Test problem}\label{s4}

In this section, we formulate a test nonlinear problem called TfC below for brevity,
which describes the minimum-time motion from a state $(x_0,y_0)$ to a state $(x_f,y_f)$
with an inequality constrained control.

The problem TfC has the following components:
\begin{itemize}
\item State vector:
 $\vec{x}=\left[\begin{array}{c}x\\y\end{array}\right]$.
Input: $\vec{u}=\left[\begin{array}{c}u\\u_d\end{array}\right]$.
\item Parameter variables: $\vec{p}=[t_f]$, where $t_f$ denotes the length
of the evaluation horizon.
\item Dynamics:  $\dot{\vec{x}}=f(\vec{x},\vec{u},\vec{p})=
\left[\begin{array}{c}(Ax+B)\cos u\\(Ax+B)\sin u\end{array}\right]$.
\item Constraints: $C(\vec{x},\vec{u},\vec{p})=[(u-c_{u})^2+u_d^2-r_{u}^2]=0$,
i.e., the control $u$ always stays within the band $c_{u}-r_{u}\leq u\leq c_{u}+r_{u}$).
\item Terminal constraints: $\psi(\vec{x},\vec{p})=\left[\begin{array}{c}
x-x_f\\y-y_f\end{array}\right]=0$ (the state should pass through the point
$(x_f,y_f)$ at $t=t_f$)
\item Objective function to minimize:
\[
J=\phi(\vec{x},\vec{p})+\int_t^{t+t_f}L(\vec{x},\vec{u},\vec{p})dt',
\]
where
\[
\phi(\vec{x},\vec{p})=t_f,\quad L(\vec{x},\vec{u},\vec{p})=-w_{d}u_d
\]
(the state should arrive at $(x_f,y_f)$ in the shortest time;
the function $L$ serves to stabilize the slack variable $u_d$)
\item Constants: $A=B=1$, $x_0=y_0=0$, $t_0=0$, $x_f=y_f=1$,
$c_{u}=0.8$, $r_{u}=0.2$, $w_{d}=0.005$.
\end{itemize}

The components of the corresponding discretized problem on the horizon are given below:
\begin{itemize}
\item the scaled horizon time $(\tau-\tau_0)/t_f\in[0,1]$ substitutes the 
original horizon time $\tau\in[\tau_0,\tau_0+t_f]$;
\item the discretized scaled horizon time is thus $\tau_i=i\Delta\tau$, where $i=0,1,\ldots,N$, and $\Delta\tau=1/N$;
\item the participating variables are the state $\left[\begin{array}{c}
x_i\\y_i\end{array}\right]$, the costate $\left[\begin{array}{c}
\lambda_{1,i}\\\lambda_{2,i}\end{array}\right]$, the control $\left[\begin{array}{c}
u_{i}\\u_{di}\end{array}\right]$, the Lagrange multipliers
$\mu_i$ and $\left[\begin{array}{c}\nu_{1}\\\nu_{2}\end{array}\right]$;
\item the state is governed by the model equation
\[
\left\{\begin{array}{l} x_{i+1}=x_i+\Delta\tau\left[p\left(Ax_{i}+B\right)\cos u_{i}\right],\\
y_{i+1}=y_i+\Delta\tau\left[p\left(Ax_{i}+B\right)\sin u_{i}\right],\end{array}\right.
\]
where $i=0,1,\ldots,N-1$;
\item the costate is determined by the backward recursion ($\lambda_{1,N}=\nu_1$,
$\lambda_{2,N}=\nu_2$)
\[
\left\{\begin{array}{l} \lambda_{1,i}=\lambda_{1,i+1}\\
\hspace{2.5em}{}+\Delta\tau\left[pA(\cos u_i \lambda_{1,i+1}+\sin u_i\lambda_{2,i+1})\right],\\
\lambda_{2,i} = \lambda_{2,i+1},\end{array}\right.
\]
where $i=N-1,N-2,\ldots,0$;
\item the equation $F(U,x_0,t_0)=0$, where
\begin{eqnarray*}
\lefteqn{U=[u_0,u_{d,0},\ldots,u_{N-1},u_{d,N-1},}\hspace*{8em}\\
&&\mu_0,\ldots,\mu_{N-1},\nu_1,\nu_2,p],
\end{eqnarray*}
has the following rows from the top to the bottom:
\[
\left\{\begin{array}{l}
\Delta\tau p\left[(Ax_i+B)\left(-\sin u_i\lambda_{1,i+1}+
\cos u_i\lambda_{2,i+1}\right)\right.\\
\hspace*{11em}\left.{}+2\left(u_i-c_{u}\right)\mu_i\right] = 0 \\
\Delta\tau p\left[2\mu_iu_{di}-w_{d}\right] = 0 \end{array}\right.
\]
\[
\left\{\;\;\Delta\tau p\left[(u_i-c_{u})^{2}+u_{di}^2-r_{u}^2\right]=0
\right.\hspace*{8em}
\]
\[
\left\{\begin{array}{l} x_N-x_r=0\\y_N-y_r=0\end{array}\right.\hspace{15em}
\]
\[
\left\{\begin{array}{l}\Delta\tau \{\sum\limits^{N-1}_{i=0}
(Ax_i+B)(\cos u_i\lambda_{1,i+1}+\sin u_i\lambda_{2,i+1})\\
\hspace{0em}{}+\mu_i[(u_i-c_u)^2+u_{di}^2-r_u^2]-w_du_{di}\}+1 = 0.\end{array}\right.
\]
\end{itemize}
Substituting  $p\mu_i$ for $\mu_i$, prior to differentiating the Lagrangian, 
leads to alternative simpler and more numerically stable, as observed in our tests, formulas, as follows
\[
\left\{\begin{array}{l}
\Delta\tau\left[p(Ax_i+B)\left(-\sin u_i\lambda_{1,i+1}+
\cos u_i\lambda_{2,i+1}\right)\right.\\
\hspace*{11em}\left.{}+2\left(u_i-c_{u}\right)\mu_i\right] = 0 \\
\Delta\tau\left[2\mu_iu_{di}-w_{d}p\right] = 0 \end{array}\right.
\]
\[
\left\{\;\;\Delta\tau\left[(u_i-c_{u})^{2}+u_{di}^2-r_{u}^2\right]=0
\right.\hspace*{8em}
\]
\[
\left\{\begin{array}{l} x_N-x_r=0\\y_N-y_r=0\end{array}\right.\hspace{15em}
\]
\[
\left\{\begin{array}{l}\Delta\tau [\sum\limits^{N-1}_{i=0}
(Ax_i+B)(\cos u_i\lambda_{1,i+1}+\sin u_i\lambda_{2,i+1})\\
\hspace{12em}{}-w_du_{di}]+1 = 0.\end{array}\right.
\]

We use the latter formulas in our numerical experiments described in the next section. 

\section{Numerical results}\label{s5}

In our numerical experiments with the TfC problem the system of linear equations
(\ref{e11}) is solved by the GMRES method. We have also tested MINRES, obtaining
the controls similar to those with GMRES, reported here. The number of evaluations of
$a(V)$ in GMRES does not exceed an a priori chosen parameter denoted by $k_{\max}$,
the error tolerance is $tol=10^{-5}$. The sampling time in the evaluation horizon is
$\Delta\tau=0.1$, the sampling time of the simulation is $\Delta t=0.02$,
and $h=10^{-5}$.

The preconditioners are constructed as follows.
At the time instances $t=jt_p$, $j=0,1,\ldots$, with an a priori chosen time increment
$t_p$ we calculate all entries of the matrix $A$ by (\ref{e10}) and its LU factorization
$A=LU$ by Gaussian elimination with partial pivoting. The computed factors $L$
and $U$ are then used in the preconditioner as follows $Tr=U^{-1}(L^{-1}r)$ for all sampling points
$t=i\Delta t$ in the interval $[jt_p,(j+1)t_p)$.
%
%The method with no preconditioning and $k_{\max}=2$ fails to compute the control within the constraints, 
%thus, the results are not shown. 

The whole set of simulations reported here consists of the following four cases:
\begin{enumerate}
 \item no preconditioning, $k_{\max}=10$;
 \item preconditioning with $t_p=0.2$ sec, $k_{\max}=1$;
 \item preconditioning with $t_p=0.4$ sec, $k_{\max}=2$;
 \item preconditioning with $t_p=0.4$ sec, $k_{\max}=10$.
\end{enumerate}

The computed results are similar in all reported cases. 
Figure  1 displays the typical CNMPC control $u$, within the constant constraints, and the time to destination $t_f$, 
both as functions of the system time in seconds, shown at the horizontal axis. 
Figure 2 shows a typical system trajectory in the $x$-$y$ plane.  

\begin{figure*}
\centering
\subfloat[NMPC control $u$ and time to destination $t_f$ for TfC
(reaches the target at $t=0.96$)]{
\label{fig1a}
\resizebox{.45\textwidth}{!}{{\scalefont{1}\includegraphics{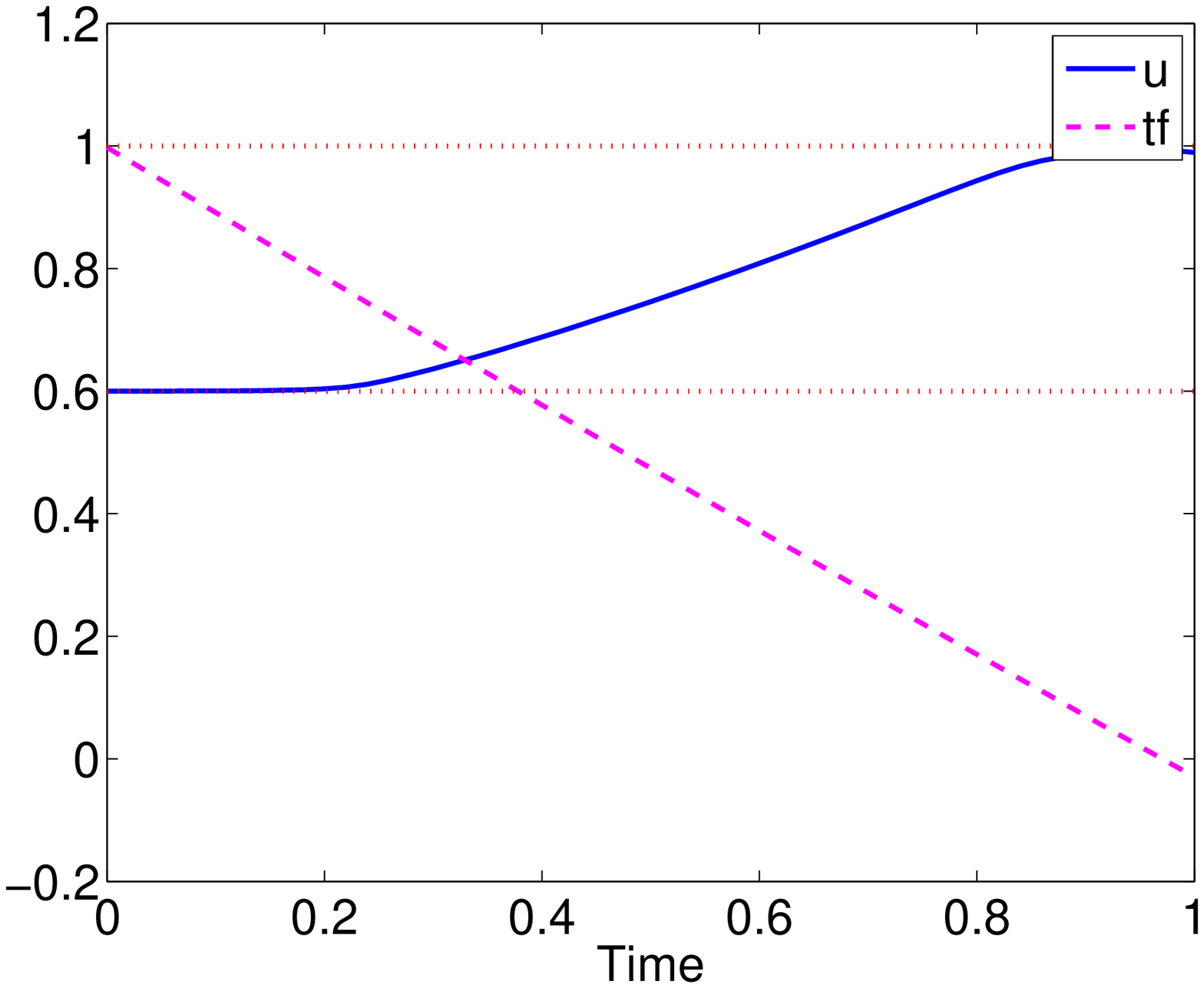}}}
}\hfill
\subfloat[TfC trajectory by NPMC]{
\label{fig1b}
\resizebox{.45\textwidth}{!}{{\scalefont{1}\includegraphics{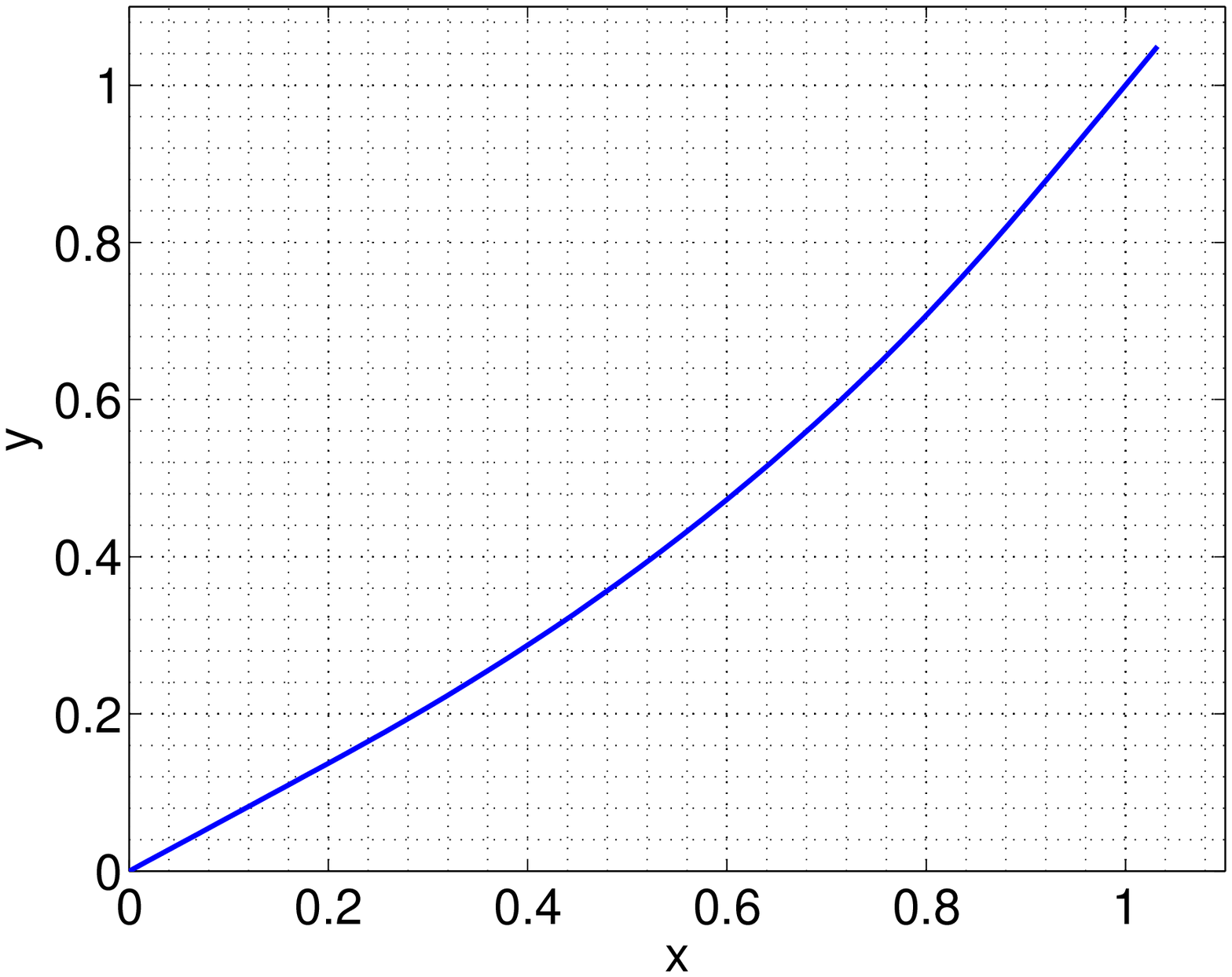}}}
}
\end{figure*} 

\begin{figure*}
\setcounter{subfigure}{2}
\centering
\subfloat[GMRES without preconditioning, $k_{\max}=10$]{
\label{fig2a}
\resizebox{.45\textwidth}{!}{{\scalefont{1}\includegraphics{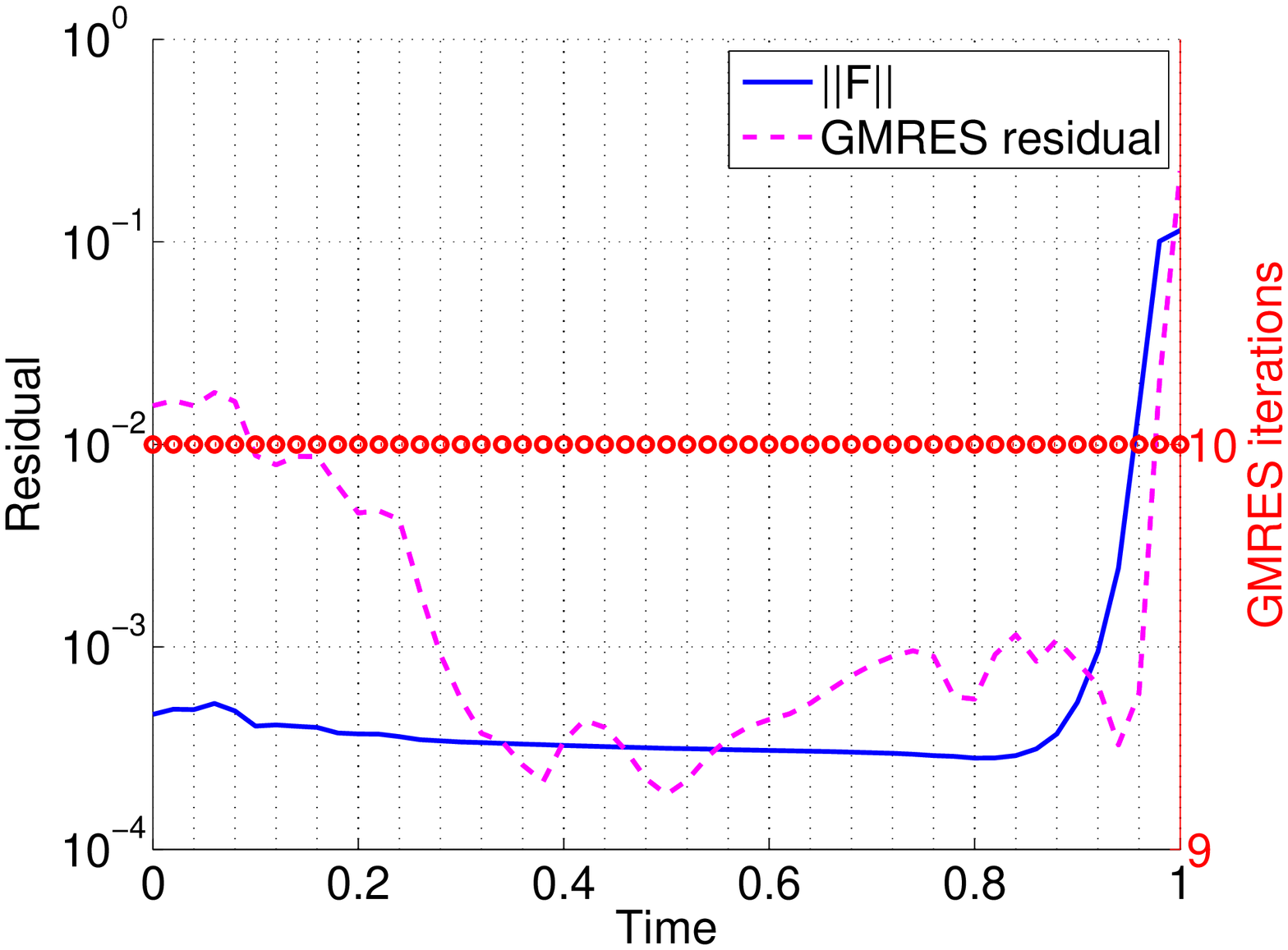}}}
}\hfill
\subfloat[GMRES with preconditionining, $t_p=0.2$ sec, $k_{\max}=1$]{
\label{fig2b}
\resizebox{.45\textwidth}{!}{{\scalefont{1}\includegraphics{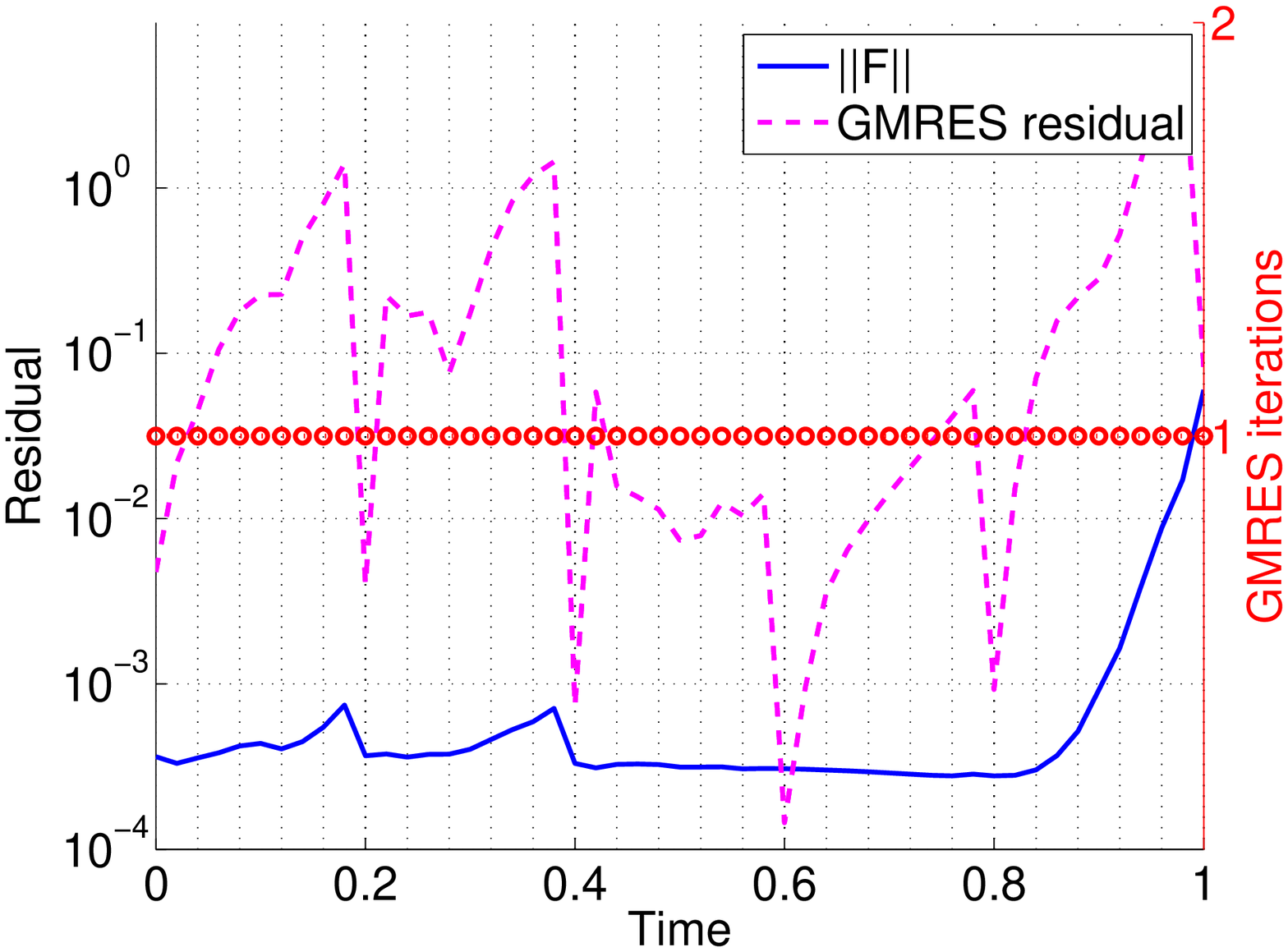}}}
}
\end{figure*}

\begin{figure*}
\setcounter{subfigure}{4}
\centering
\subfloat[GMRES with preconditionining, $t_p=0.4$ sec, $k_{\max}=2$]{
\label{fig3a}
\resizebox{.45\textwidth}{!}{{\scalefont{1}\includegraphics{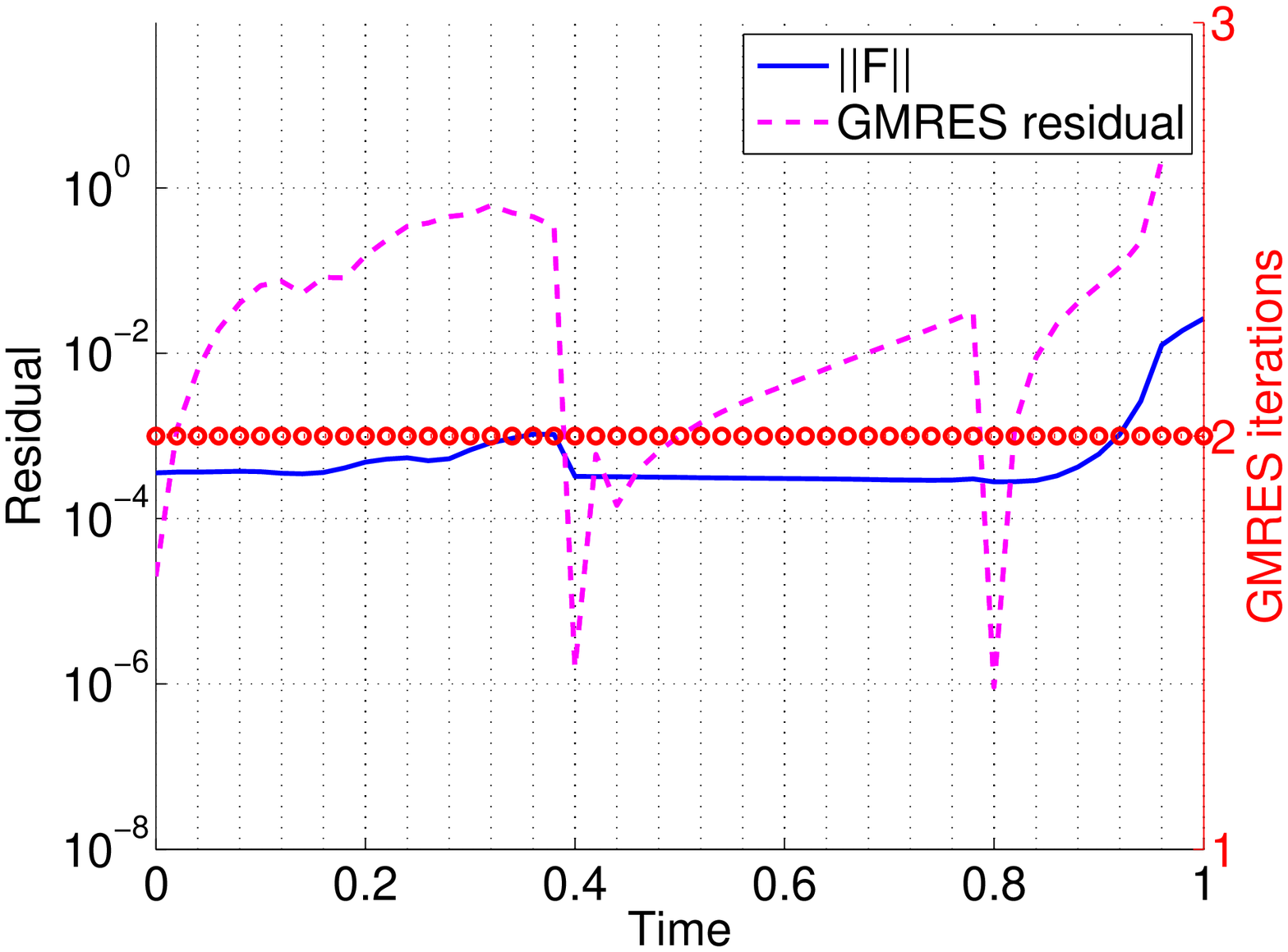}}}
}\hfill
\subfloat[GMRES with preconditioning, $t_p=0.4$ sec, $k_{\max}=10$]{
\label{fig3b}
\resizebox{.45\textwidth}{!}{{\scalefont{1}\includegraphics{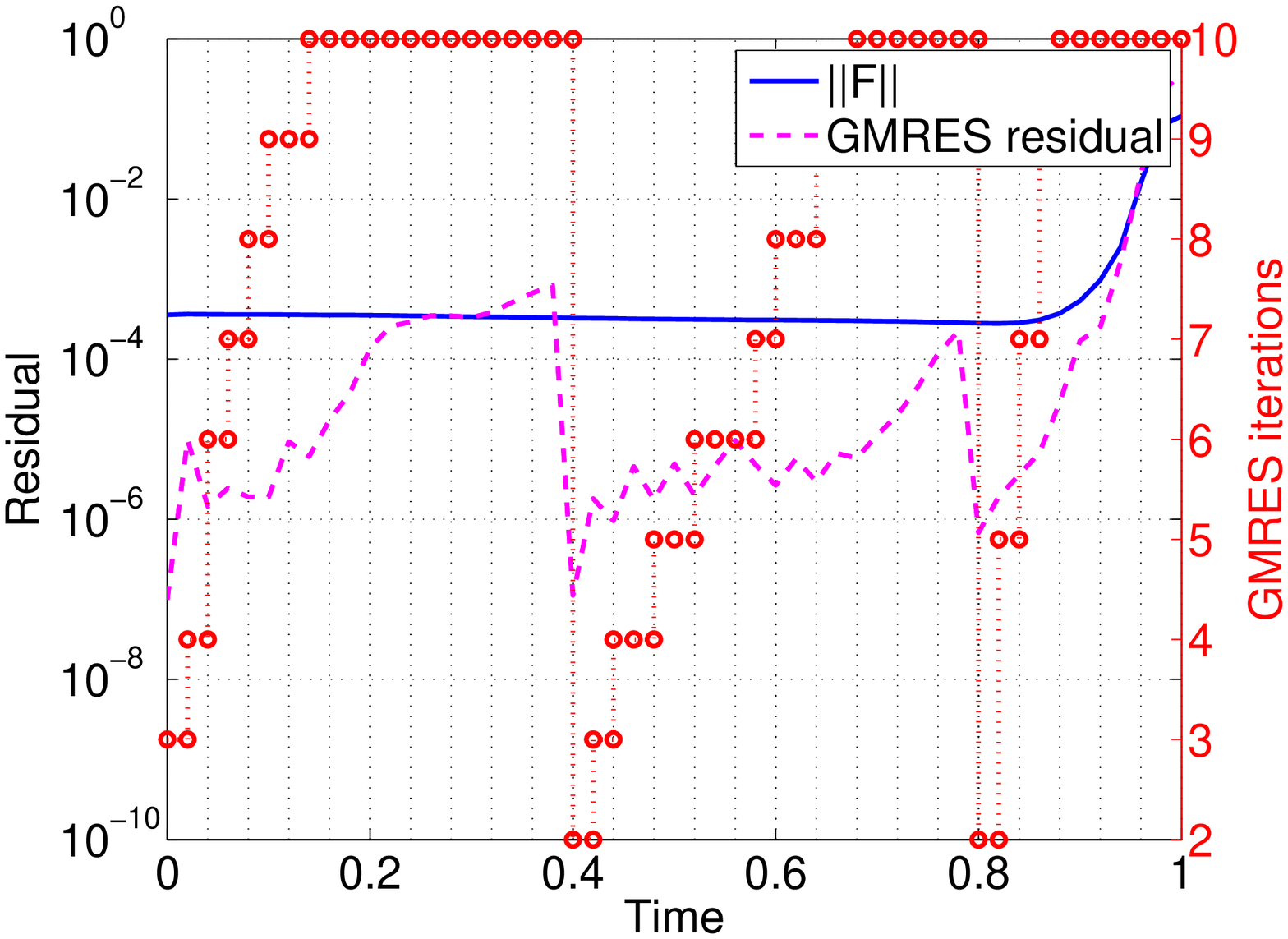}}}
}
\end{figure*}

Figures 3--6 show the value of $\|F\|$, which we want to be vanished, and the GMRES residual (the left vertical axis) and 
the number of the actually performed GMRES iterations  (the right vertical axis)
at every system time step for all four cases, where the horizontal axis represents the system time in seconds.
Figure 3 corresponds to the GMRES iterations without preconditioning. Figures 4-6 
involve the preconditioner, recalculated with various frequencies, determined by the time increment
$t_p$, and for different $k_{\max}$
ranging from $1$ to $10$.

In Figure 3, the number of the actually performed GMRES iterations without preconditioning is always the maximum allowed in this test  $k_{\max}=10$. We use this test as a baseline for comparisons. 

We first point out a good behavior of the preconditioned GMRES even with $k_{\max}=1$ and where 
the preconditioner is reconstructed once each $t_p=0.2$ sec, see Figure 4.
This clearly demonstrates the fact that preconditioning reduces
the number of evaluations of the vector function $F(U,x,t)$.

The effect of increasing the maximum number $k_{\max}$ of GMRES steps
is seen by comparing Figures 4-6. 
Specifically, in Figure 4, $t_p=0.2$ sec and $k_{\max}=1$, compared to 
$t_p=0.4$ sec and $k_{\max}=2$ in Figure~5, i.e., we can recompute 
the preconditioner twice less frequently at the cost of increasing  $k_{\max}$
from $1$ to $2$, and we observe a slightly better quality of the solution, as measured by  
the generally smaller values of $\|F\|$ and the GMRES residual (the left vertical axis).

In Figure 6, the preconditioner is recomputed as frequent as in Figure 5, but the largest allowed 
number of GMRES iterations is increased from $k_{\max}=2$ to $k_{\max}=10.$ 
We observe in Figure 6 that GMRES often activates the default tolerance stopping criteria 
for the residual norm smaller than $10^{-5}$, before maxing out the allowed number of iterations $k_{\max}$. 
Overall, this leads to a generally much smaller residual in Figure 6 compared to that in Figure 5. 
However, the most decisive quantity $\|F\|$ behaves similar both in Figures 5 and 6, and the 
computed controls are so similar that the increase of $k_{\max}$ from $2$ to $10$
may be unnecessary.  

Efficiency of preconditioning is illustrated by comparing
Figures 3 and 5, where the number of iterations is reduced five times 
giving similar/smaller values of   $\|F\|$. 

In minimum-time optimal control problems, the length of
the evaluation horizon gets smaller as the state $(x,y)$ approaches
the goal position. Near the goal position $(1,1)$
the control has less capability (controllability)
to direct the state towards the goal because of short time for control.
This makes the equation $F(U)=0$ more difficult
for numerical solution, thus, $\|F\|$ increases near the goal position, as 
seen in Figures~3--6.

\section*{Conclusions}

Time-optimal problems are practically important, giving optimal solutions 
    for guidance, navigation and control, 
    which can be used for vehicles, trains, etc.
Due to heavily nonlinear equations and highly coupled variables, the 
time-optimal problems are difficult to solve numerically. 
We present an apparently first successful extension of CNMPC  
for real-time control of such problems. 
Our numerical experiments demonstrate dramatic acceleration of convergence
of iterations without sacrificing control quality, if proper preconditioning
is used. The proposed concurrent construction of the preconditioner can be
trivially efficiently implemented in parallel on controllers having multiple
processing units, such as multi-core, graphics processing units, and  modern
field-programmable gate arrays. Replacing GMRES with the MINRES iterative solver
may help reducing controller memory requirements and increasing the speed of convergence. 
Our algorithm, including the preconditioner setup implemented
in parallel and the iterative solver, can significantly speed up the calculation
of the control, compared to traditional sequential CNMPC algorithms, thus
allowing to control system with faster dynamics. 
Our future work concerns analyzing MINRES, as a possible 
replacement of GMRES, and developing efficient preconditioners, with faster on-line setup and application,
within the framework of CNMPC. 
%\pagebreak

\end{document}